\documentclass{agtart_a}
\pdfoutput=1
\usepackage[matrix,arrow,curve]{xy}

\title{Dyer--Lashof--Cohen operations in Hochschild cohomology}

\author{Victor Tourtchine}
\givenname{Victor}
\surname{Tourtchine}
\address{Institut de Math\'ematique Pure et Appliqu\'ee\\
Universit\'e Catholique de Louvain\\\newline
Chemin du Cyclotron 2\\
B-1348 Louvain-la-Neuve\\ 
Belgium}
\email{vitia-t@yandex.ru}
\urladdr{}

\volumenumber{6}
\issuenumber{}
\publicationyear{2006}
\papernumber{33}
\startpage{875}
\endpage{894}

\doi{}
\MR{}
\Zbl{}


\received{26 April 2005}
\revised{20 November 2005}
\accepted{5 April 2006}
\published{24 July 2006}
\publishedonline{24 July 2006}
\proposed{}
\seconded{}
\corresponding{}
\editor{}
\version{}

\arxivreference{math.RA/0504017}

\let\xysavmatrix\xymatrix
\def\xymatrix{\disablesubscriptcorrection\xysavmatrix}
\AtBeginDocument{\let\bar\wbar\let\tilde\wtilde\let\hat\what}

\makeop{sign}
\makeop{lf}
\makeop{deg}

\makeatletter
\def\cnewtheorem#1[#2]#3{\newtheorem{#1}{#3}[section]
\expandafter\let\csname c@#1\endcsname\c@theorem}

\newtheorem*{theorem*}{Theorem}
\newtheorem{theorem}{Theorem}[section]
\cnewtheorem{lemma}[theorem]{Lemma}
\cnewtheorem{statement}[theorem]{Statement}
\cnewtheorem{conjecture}[theorem]{Conjecture}
\cnewtheorem{proposition}[theorem]{Proposition}
\cnewtheorem{corollary}[theorem]{Corollary}

\theoremstyle{definition}
\cnewtheorem{remark}[theorem]{Remark}
\cnewtheorem{definition}[theorem]{Definition}
\cnewtheorem{example}[theorem]{Example}
\makeatother

\def\mod{\mathop{\rm mod}\nolimits}

\newcommand\N{{\mathbb N}}

\newcommand\HH{{\mathcal H}}

\newcommand\kk{\mathbb{K}}

\newcommand\Assoc{{\mathcal {ASSOC}}}
\newcommand\Comm{{\mathcal {COMM}}}

\newcommand\Lie{{\mathcal {LIE}}}

\newcommand\Gerst{{\mathcal {GERST}}}
\renewcommand\O{{\mathcal O}}
\newcommand\Brace{{\mathcal {BRACE}}}
\newcommand\prelie{{\mathcal {PL}}}
\newcommand\G{G_\infty}
\newcommand\vertical{{\mathcal {VERT}}}

\newcommand\FF{{\mathcal {F}}}
\newcommand\CC{{\mathcal {C}}}
\newcommand\BB{{\mathcal {B}}}
\newcommand\X{{\mathcal {X}}}
\newcommand\SSS{{\mathcal {S}}}

\def\ad{{\rm ad}}

\numberwithin{equation}{section}

\renewcommand{\thefootnote}{\fnsymbol{footnote}}

\begin{document}

\begin{asciiabstract}
We give explicit formulae for operations in Hochschild cohomology which
are analogous to the operations in the homology of double loop spaces.
As a corollary we obtain that any brace algebra in finite characteristic
is always a restricted Lie algebra.
\end{asciiabstract}

\begin{abstract}
We give explicit formulae for operations in Hochschild cohomology which
are analogous to the operations in the homology of double loop spaces.
As a corollary we obtain that any brace algebra in finite characteristic
is always a restricted Lie algebra.
\end{abstract}

\maketitle

\setcounter{footnote}{0}
\renewcommand{\thefootnote}{\arabic{footnote}}
\section{Introduction}\label{s0}

\subsection{Homology operations}\label{s01}
It is already well known that the Hochschild cochain complex of
an associative algebra can be endowed with an action of an operad
quasi-isomorphic to the chain operad of little squares. This statement
is called ``Deligne's conjecture''. Over $\Z$ this result is due to
J\,E~McClure and J\,H~Smith~\cite{McSm1,McSm2}, and  to M~Kontsevich
and Ya~Soibelman~\cite{KS}, see also the proof of C~Berger and
B~Fresse~\cite{BF}, which appeared later. In characteristic zero there
are several proofs; see for example M~Kontsevich \cite{K2}, D\,E~Tamarkin
\cite{Tam1,Tam2} and A\,A Voronov \cite{Vor}.  This result implies that
in Hochschild cohomology one can define the same homological operations
as for double loop spaces.  Homological operations for the iterated loop
spaces are well known~\cite{Co}: In the case of double loops, one has
a Pontryagin multiplication, a Browder operator (degree one bracket),
and also two non-trivial Dyer--Lashof operations (following F~Cohen we
denote them by $\xi_{1}$ and $\zeta_{1}$):

Over $\Z_{2}$:
 \begin{equation}
 \xi_{1}\co H_{k}(\Omega^2X,\Z_{2})\to H_{2k+1}(\Omega^2X,\Z_{2}).
 \label{eq01}
 \end{equation}
 Over $\Z_{p}$, $p$ being any odd prime:
 \begin{gather}
 \xi_{1}\co H_{2k-1}(\Omega^2X,\Z_{p})\to H_{2pk-1}(\Omega^2X,\Z_{p}),
 \label{eq02}
 \\
 \zeta_{1}\co H_{2k-1}(\Omega^2X,\Z_{p})\to H_{2pk-2}(\Omega^2X,\Z_{p}).
 \label{eq03}
 \end{gather}
Over $\Z_{2}$, operation $\xi_{1}$ was introduced by  S~Araki and T~Kudo~\cite{AK2}. Over $\Z_{p}$, operations $\xi_{1}$ and $\zeta_{1}$
were introduced by F~Cohen~\cite{Co}. All the other operations are some compositions of the 
above~\cite{Co}.

The study of $\Z_{p}$--homology operations, $p>3$, for iterated $d$--loop
spaces was initiated by E~Dyer and R\,K~Lashof~\cite{DyLash}. But they
found only  part of the homological
operations. For example, in our case $d=2$ their method did not recover operations~\eqref{eq02}, \eqref{eq03}. A complete list of operations together
with all the relations was given by F~Cohen, cf \cite{Co}.

The above homological operations appear via an action of the operad $\CC_{2}$ of little squares, and correspond to specific cycles of equivariant homology of 
$\CC_{2}$. To be precise, choose a homology class $\alpha\in H_{*}(\CC_{2}(n)/S_{k_1}\times\cdots\times S_{k_\ell},W)$, where $k_{1}+\cdots +k_{\ell}=n$,
\begin{equation}
W=\otimes_{i=1}^{\ell}(\sign_{i})^{\otimes d_i},
\label{eq04}
\end{equation}
with each factor $\sign_i$ being a sign representation of $S_{k_{i}}$. To this cycle we can associate a homological
operation
$$
\alpha\co H_{d_{1}}(\Omega^2X)\times H_{d_{2}}(\Omega^2X)\times\cdots\times H_{d_{\ell}}(\Omega^2X)\to
H_{k_{1}d_{1}+\cdots+k_{\ell}d_{\ell}+\deg(\alpha)}(\Omega^2X),
$$
which we  denote by the same letter $\alpha$. If $k_{1}=k_{2}=\cdots =k_\ell=1$, then $\alpha$ is a multilinear operation.

The Pontryagin product $*$, the Browder operator $[\, .\, , \, . \, ]$, $\xi_{1}$ and $\zeta_{1}$  correspond to cycles
$$
\begin{array}{ccl}
* &\in&H_{0}(\CC_{2}(2))=H_{0}(S^1);\\
{[\, .\, , \, . \, ]}&\in&H_{1}(\CC_{2}(2))=H_{1}(S^1);\\
\xi_{1}&\in&H_{p-1}(\CC_{2}(p)/S_{p},\pm\Z_{p});\\
\zeta_{1}&\in&H_{p-2}(\CC_{2}(p)/S_{p},\pm\Z_{p}),\,  p>2.
\end{array}
$$
For Hochschild complexes, operations $*$, $[\, .\, , \, . \, ]$ are respectively the cup-product and the Gerstenhaber bracket~\cite{Ge}. The aim of this paper
is to give explicit formulae for operations $\xi_{1}$ and $\zeta_{1}$. 
The author used them in~\cite{T5} to describe 
the Hochschild homology of the Poisson algebras operad and of the Gerstenhaber algebras operad in the bigradings spanned by the operad of Lie algebras.

\subsection{Results}\label{s02}
The results of the paper are given by Theorems~\ref{t32}, \ref{t44}, \ref{t75}, \ref{t78}.  \fullref{t75} is also an obvious consequence of a result of  B~Fresse,
see \fullref{r64}.

\subsection{Notations}\label{s02'}

 The letter $p$ always denotes  a prime number.

%\smallskip

\noindent We suppose that  the operads and  the homology of spaces are defined over some commutative ring $\kk$, which is sometimes $\Z_{p}$. 

%\smallskip

By $H^{\lf}_{*}(\, .\, , L)$,  $H_{\lf}^{*}(\, .\, , L)$ we denote locally finite singular (co)homology with coefficients in a local system~$L$.

By $\wwbar X$ we denote the one-point compactification of a space $X$.

\subsection{Acknowledgement}\label{s03}
The author is grateful to the Universit\'e Catholique de Louvain for
the hospitality and for its friendly atmosphere.  The author thanks
P~Lambrechts, L~Schwartz, Y~Felix, F~Cohen, A~Voronov, J\,E~McClure,
D~Chataur, S~Kallel, G~Sharygin, M~Kontsevich, M~Livernet, and V~Ostrik
for interesting discussions and communications.

The author is also grateful to the referee for useful suggestions and
comments and to Hal Sadofsky for correcting the English.

While accomplishing this work the author was partially supported by the
grants NSH-1972.2003.01, RFBR~05-01-01012a.

\section{Operad $\SSS_{2}$}\label{s2}
There is a natural differential graded operad acting on Hochschild
cohomology complexes. This operad was denoted $\G$ by Gerstenhaber and
Voronov~\cite{GV},
$\HH$ by McClure and Smith~\cite{McSm1},
and $F_{2}\X$ by Berger and Fresse~\cite{BF}. We adopt the notation of
McClure and Smith~\cite{McSm2} and denote this operad by $\SSS_{2}$. This operad is generated by brace operations
\begin{equation}
\{\,\}_{n+1}\co x_{1}\otimes x_{2}\otimes\dots\otimes x_{n+1}\mapsto x_{1}\{x_{2},\dots,x_{n+1}\}, \quad n\geq 1,
\label{eq21}
\end{equation}
augmenting the degree by $n$, and by an associative cup product $*$. The relations between these operations are standard brace relations, associativity
of $*$, and standard relations between $*$ and $\{\,\}_{n}$. We refer the reader to one of the above papers for an explicit description of
this operad together with the differential on it. We suppose that the degree zero component $\SSS_{2}(0)$ of this operad is trivial. The minimal degree 
part of this operad is the associative algebras operad: $(\SSS_{2}(n))_{0}=\Assoc(n)$. The maximal degree part is the operad of shifted brace algebras:
$(\SSS_{2}(n))_{n-1}=\Brace_{1}(n)$, ie the operad of brace algebras with operations $\{\,\}_{n}$, $n\geq2$, of degree $n-1$. We have  in particular that
$(\SSS_2(n))_*$ vanishes in degree $*<0$ and in degree $*>n-1$.

\begin{theorem*}[McClure and Smith \cite{McSm1}]
The  operad $\SSS_{2}$ is quasi-isomorphic to the operad  $S_{*}(\CC_{2})$ of singular chains of little squares.
\end{theorem*}

It can be easily seen that on each component $\SSS_2(n)$  the action of the symmetric group 
$S_{n}$ is free. 

The following lemma is a particular case of \cite[Theorem 10.4.8]{Weibel}.

\begin{lemma}\label{l20}
If  left-bounded complexes of projective $G$--modules ($G$ being a finite group) are 
quasi-isomorphic then they are homotopy equivalent.\footnote{By ``quasi-isomorphic'' and ``homotopy equivalent'' we understand quasi-isomorphic and homotopy equivalent in the category of complexes of 
$G$--modules.}
\end{lemma}

As a consequence of the result of McClure--Smith and of \fullref{l20} we get:

\begin{corollary}\label{c21}
One has a natural isomorphism
$$
H_{*}(\SSS_{2}(n)\otimes_{S_{k_1}\times\cdots\times S_{k_\ell}}W)\simeq H_{*}(\CC_{2}(n)/S_{k_1}\times\cdots\times S_{k_\ell},W)
$$
for any representation $W$ of $S_{k_1}\times\cdots\times S_{k_\ell}$. 
\end{corollary}

We will need this corollary only in the case $W$ is of type~\eqref{eq04}. 

\section{Explicit formulae for $\xi_{1}$ and $\zeta_{1}$}\label{s3}
For any element $x$ of a brace algebra, denote by $x^{[n]}$ the following expression:
\begin{equation}
x^{[k]}:=x\underbrace{\{ x\}\dots\{ x\} }_{\text{$k-1$ times}}. 
\label{eq31}
\end{equation}
For example,
$$
x^{[3]}=x\{x\}\{x\}=
\begin{cases}
  x\{x\{x\}\}+2x\{x,x\}, & \text{if $\deg(x)-\deg\{\}_2$ is even};\\
  x\{x\{x\}\}, & \text{if $\deg(x)-\deg\{\}_2$ is odd.}
  \end{cases}
$$

\begin{theorem}\label{t32}
The following operations are the Dyer--Lashof--Cohen operations induced by the action of the operad $\SSS_{2}$:
\begin{align}
\xi_{1}(x)&=x^{[p]}, &&\text{$p(\deg(x)-1)$ being even;} \label{eq33}\\
\zeta_{1}(x)&=\sum_{i=1}^{p-1}\frac{(-1)^{i}}{i}x^{[i]}* x^{[p-i]},
&&\text{$p\cdot \deg(x)$ being odd}.\label{eq34}
\end{align}
\end{theorem}
 
\begin{remark}\label{r32}
For $p=2$, $\xi_1(x)=x\{x\}$.  This result is due to C~Westerland~\cite{Craig}.
\end{remark}
 
\begin{example} 
For $p=3$,
\begin{align*}
\xi_1(x)&=x\{x\}\{x\}=x\{x\{x\}\}+2x\{x,x\},\\
\zeta_1(x)&=-x*x\{x\}-x\{x\}*x.
\end{align*}
\end{example}
 
\begin{proof}[Proof of \fullref{t32}] We proved in~\cite[Section 11]{T5}, that if $x$ is an odd degree cycle in a Hochschild complex (or of any degree with characteristic $p=2$), then the following formula holds:\footnote{This formula can be proved by induction over $n$, it is also a consequence of
 equality~\eqref{eq47}, which we prove in \fullref{s4}.}
\begin{equation}
\partial(x^{[n]})=-\sum_{i=1}^{n-1} {n\choose i }x^{[i]}*x^{[n-i]}.
\label{eq34'}
\end{equation}
It follows from~\eqref{eq34'} that~\eqref{eq33} defines some
  cycle $\alpha$ of the complex $\SSS_{2}(p)\otimes_{S_p}(\pm\Z_{p})$.   This cycle defines a non-trivial homology class, since it is in the maximal 
  degree $*=p-1$ of the above complex. It is well known (see Cohen
  \cite{Co}, Vassiliev \cite{V0} and Markaryan \cite{Ma}) that
  $$
  H_{*}(B(p,\R^2),\pm\Z_{p})=H_{*}(\CC_{2}(p)/S_{p},\pm\Z_{p})=
  \begin{cases}
  \Z_{p},&\text{$*=p-2$ or $p-1$,}\\
  0,&\text{otherwise,}
  \end{cases}
  $$
where $B(p,\R^2)$ denotes as usual the configuration space of cardinality $p$ subsets of $\R^2$.

By \fullref{c21},
$$
H_{*}(\SSS_{2}(p)\otimes_{S_p}(\pm\Z_{p}))=\begin{cases}
  \Z_{p},&\text{$*=p-2$ or $p-1$,}\\
  0,&\text{otherwise.}
  \end{cases}
$$
It means that operation~\eqref{eq33} is a multiple of $\xi_{1}$,
$$
\alpha=\lambda\xi_{1},
$$
for some coefficient $\lambda\neq 0$. 
%So, the equality~\eqref{eq76} must hold with coefficient $\lambda$, see \fullref{s6}. But it can be easily seen
%that this coefficient can be only one. It implies the result for operation~\eqref{eq33}.
We will prove in \fullref{s6} that this coefficient $\lambda$ is exactly one.

To see that formula~\eqref{eq34} defines operation $\zeta_{1}$, it is sufficient to show that cycle~\eqref{eq34} is the image of the Bockstein homomorphism
$\beta$ of the cycle~\eqref{eq33}.\footnote{We warn the reader that  the cycle $\zeta_1x$ in Hochschild cohomology is not simply $\beta(\xi_1x)$, but is related to the above
by the following formula of F~Cohen~\cite{Co}:
$$\zeta_{1}x=\beta(\xi_{1}x)- (\ad^{p-1}x)(\beta x).$$
Operator $\ad\, x$ is the adjoint action $[x,\, .\, ]$. Surprisingly, the above formula holds on the level of chains, see \fullref{p65}.} This follows from~\eqref{eq34'}, and also from the equality
$$
-\frac{(p-1)!}{i!(p-i)!}\equiv\frac{(-1)^i}{i}\mod p.\proved
$$
\end{proof}

\section{Quasi-isomorphism $\FF_{2}\to\SSS_{2}$}\label{s4}

In this  section we will give another proof of \fullref{t32}. This construction is interesting in itself.

It turns out that complexes $\SSS_{2}(n),\, n\geq 1$, are too big, and that they contain much smaller 
quasi-isomorphic  subcomplexes $\FF_{2}(n),\, n\geq 1$.
Complexes $\FF_{2}(n),\, n\geq 1$, do not form an operad, but they are freely acted on by $S_{n}$, and so \fullref{l20} can be applied. These
complexes have a geometric origin.

We will first define complexes $\BB_{2}(n,\kk)$, $\BB_{2}(n,\pm\kk)$, which are in fact \mbox{$\FF_{2}(n)\otimes_{S_{n}}\kk$},
$\FF_{2}(n)\otimes_{S_{n}}\pm\kk$.

Consider the space $B(n,\R^2)$ of cardinality $n$  subsets of $\R^2$. The pace $B(n,\R^2)$ is homotopy equivalent to $\CC_{2}(n)/S_{n}$.
By Pontryagin duality,
\begin{gather}
H_{*}(B(n,\R^2),\kk)\simeq  H_{\lf}^{2n-*}(B(n,\R^2),\kk);\label{eq41}\\
H_{*}(B(n,\R^2),\pm\kk)\simeq H_{\lf}^{2n-*}(B(n,\R^2),\pm\kk);\label{eq42}
\end{gather}
where $H_{\lf}^{*}(\, .\, , L)$ denotes locally finite singular cohomology with coefficients in a local system~$L$.

To compute the right hand side of~\eqref{eq41},~\eqref{eq42}, one can use the following cellular decomposition of the one point compactification 
$\overline{B(n,\R^2)}$,~\cite{F,Vain,V0,Ma}.  Let  $A=\{a_1,a_2,\dots,a_n\}$ be a point of $B(n,\R^2)$. We will assign to $A$ its {\it
index}. This is a system of numbers $(k_1,k_2,\dots,k_\ell)$ satisfying $k_1+k_2+\dots +k_\ell=n$, where $k_1$ is the number of
elements of $A$ with the minimal value of the  first coordinate $x$; $k_2$ is the number of elements of
$A$ with next value of $x$, and so on$\dots$ Points with the same index $(k_1,k_2,\dots,k_\ell)$ form a cell,
that we denote by $e(k_1,k_2,\dots,k_\ell)$. 

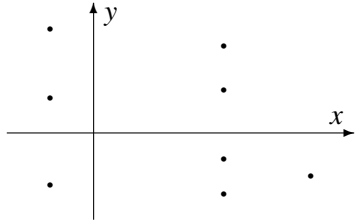
\begin{figure}[h!t]
\unitlength=0.3em
\begin{center}
\begin{picture}(40,25)

\put(0,10){\vector(1,0){40}}
\put(10,0){\vector(0,1){25}}
\put(5,4){\circle*{0.5}}
\put(5,14){\circle*{0.5}}
\put(5,22){\circle*{0.5}}

\put(25,15){\circle*{0.5}}
\put(25,20){\circle*{0.5}}
\put(25,7){\circle*{0.5}}
\put(25,3){\circle*{0.5}}

\put(35,5){\circle*{0.5}}
\put(37,11){$x$}
\put(11,23){$y$}
\end{picture}
\end{center} 
\caption{Point of the cell $e(3,4,1)$ of $\overline{B(8,\R^2)}$.}  
\end{figure}

All such cells together with the infinite point provide a cell
decomposition of $\overline{B(k,\C)}$. These cells bound to each other by the rule
\begin{multline}
\partial e(k_1,\dots,k_\ell)=\\
\sum_{i=1}^{\ell-1}(-1)^{i-1}{{k_i{+}k_{i+1}}\choose{k_i}}
e(k_1,\dots,k_{i-1},k_i{+}k_{i+1},k_{i+2},\dots,k_\ell)
\label{eq42'}
\end{multline}
for twisted coefficients $\pm\kk$; and by the rule 
\begin{multline}
\partial  e(k_1,\dots,k_\ell) =\\
\sum_{i=1}^{\ell-1}(-1)^{s_{i}}{{k_i{+}k_{i+1}}\choose{k_i}}_{-1}
e(k_1,\dots,k_{i-1},k_i{+}k_{i+1},k_{i+2},\dots,k_\ell)
\label{eq42''}
\end{multline}
for constant coefficients, where $s_{i}=i-1+k_1+k_2+\dots+k_i$,
$$
{{k+\ell}\choose k}_{-1}=
\begin{cases}
0,&\text{$k$ and $\ell$ are odd;}\\
\left({{\left[\frac{k+\ell}2\right]}\atop{\left[ \frac k2\right]}}\right), & \text{otherwise}.
\end{cases}
$$
Since we are interested not in the homology but in the cohomology
$H_{\lf}^*(B(n,\R^2),-)$, we need to consider the \emph{duals} of
the above complexes.  Let us denote these duals by $\BB_{2}(n,\pm\kk)$,
and $\BB_{2}(n,\kk)$ respectively.  The elements of the dual basis in
these dual complexes will be denoted by $\epsilon(k_1,\ldots,k_\ell)$.

Now, define complexes $\FF_{2}(n)$.
Consider the space $F(n,\R^2)$ of $n$ distinct points in $\C$. Obviously, 
$F(n,\R^2)/S_{n}=B(n,\C)$. By Poincar\'e duality,
\begin{equation}
H_{*}(F(n,\R^2))\simeq H_{\lf}^{2n-*}(F(n,\R^2)).
\label{eq42'''}
\end{equation}
We will consider a cell decomposition of $\overline{F(n,\C)}$, which is a preimage of the above cell decomposition
of $\overline{B(n,\C)}$. Explicitly, each cell $e(\sigma;k_{1},\dots,k_{\ell})$ of $F(n,\C)$ is encoded by a permutation $\sigma\in S_{n}$ and
a sequence $(k_{1},\dots,k_{\ell})$ of positive integers, such that $k_{1}+\dots +k_{\ell}=n$. A point $\bar A=(a_{1},a_{2},\dots,a_{n})\in F(n,\R^2)$
belongs to $e(\sigma;k_{1},\dots,k_{\ell})$ if $A=\{a_{1},a_{2},\dots,a_{n}\}\in B(n,\R^2)$ belongs to $e(k_{1},\dots,k_{\ell})$, and the order 
of indices is $\sigma_1,\sigma_2,\dots,\sigma_n$ when the points $a_{1},\dots,a_{n}$ are lexicographically ordered.

\begin{figure}[!ht]
\unitlength=0.3em
\begin{center}
\begin{picture}(20,20)

\put(0,10){\vector(1,0){20}}
\put(10,0){\vector(0,1){20}}
\put(5,7){\circle*{0.5}}
\put(5.4,6.7){$3$}

\put(13,5){\circle*{0.5}}
\put(13.4,4.7){$1$}
\put(13,13){\circle*{0.5}}
\put(13.4,12.7){$2$}
\put(19,11){$x$}
\put(11,18){$y$}
\end{picture}
\end{center} 
\caption{Point of the cell $e(\sigma;1,2)$, where $\sigma=(3,1,2)$.}  
\end{figure}
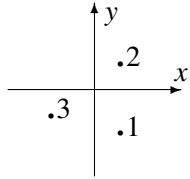

The differential of a cell is
$$\partial e(\sigma;k_{1},\dots,k_{\ell})=
\sum_{i=1}^{\ell-1}(-1)^{k_{1}+\dots +k_{i}}\hspace{-12pt}\sum_{\sigma'\in
S(i,i+1)}\hspace{-10pt}(-1)^{|\sigma'|}
e(\sigma'\sigma;k_{1},\dots,k_{i}+k_{i+1},\dots,k_{\ell}),$$
where $S(i,i+1)$ is a subset of $S_{n}$ of all shuffles of $k_{i}$ consecutive points starting from $k_{1}+\dots +k_{i-1}+1$, with 
$k_{i+1}$ consecutive points starting from $k_{1}+\dots +k_{i}+1$.

We denote by $\FF_{2}(n)$ the \emph{dual} of the above complex. Its  basis elements will be denoted by  $\epsilon(\sigma;k_{1},\dots,k_{\ell})$.

Define an inclusion
$$
\xymatrix{
I\co\FF_{2}(n)\ar@{^{(}->}[r]&\SSS_{2}(n),
}
$$
by the formula
\begin{multline*}
I(\epsilon(\sigma;k_{1},\dots,k_{\ell}))=\\
(-1)^{n-\ell}x_{\sigma_1}\{x_{\sigma_2}\}\{x_{\sigma_3}\}\ldots\{x_{\sigma_{k_{1}}}\}*
 x_{\sigma_{k_{1}+1}}\{x_{\sigma_{k_{1}+2}}\}\ldots\{x_{\sigma_{k_{1}+k_{2}}}\}*\cdots\\
 \cdots
 *x_{\sigma_{k_{1}+\cdots+k_{\ell-1}+1}}\{x_{\sigma_{k_{1}+\cdots+k_{\ell-1}+2}}\}\ldots\{x_{\sigma_n}\}.
\end{multline*}

\begin{lemma}\label{l45}
The map $I$ is a morphism of complexes.  
\end{lemma}

\begin{proof}[Proof of \fullref{l45}] The differential of the cup-product $*\in\SSS_{2}(2)$ is zero. It implies
that it is enough to prove the identity\footnote{This formula, applied to the case
$x_1=x_2=\cdots=x_n=x$ of odd degree (this affects signs), implies identity~\eqref{eq34'}.}
\begin{multline}
\partial x_{1}\{x_2\}\{x_{3}\}\dots\{x_{n}\}=\\
\sum_{ {I\sqcup J=\{1\ldots n\}\,\, I\neq\emptyset\neq J }\atop { {I=\{i_{1}<i_{2}<\dots <i_{k}\} } \atop {J=\{j_{1}<j_{2}<\dots <j_{n-k}\}}  }  }(-1)^S
x_{i_1}\{x_{i_2}\}\{x_{i_3}\}...\{x_{i_k}\}* x_{j_1}\{x_{j_2}\}\{x_{j_3}\}...\{x_{j_{n-k}}\}.\label{eq47}
\end{multline} 
Here $S=|\sigma(I,J)|+k+1$ with $|\sigma(I,J)|$ being the sign of the corresponding shuffle permutation. The sum is taken over all possible partitions of $\{1\ldots n\}$ into two non-empty subsets $I$ and $J$. Identity~\eqref{eq47} can be proven by induction over $n$. For $n=1$, it is evident. 
The inductive step
follows from the identities
\begin{align}
\partial(A\{ B\})&=(\partial A)\{ B\}-(-1)^{\deg(A)-1}A\{ \partial B\}
  \label{eq48}\\
&\qquad\qquad+(-1)^{\deg(A)}(A* B-(-1)^{\deg(A)\deg(B)}B* A),\notag\\
(A* B)\{ C\}&=A*(B\{ C\})+(-1)^{\deg(B)(\deg(C)-1)}(A\{ C\})* B.
\label{eq49}
\end{align}
For example,
\begin{equation}
\begin{aligned}
\partial (x_{1}\{x_2\}\{x_{3}\})
&=\partial(x_{1}\{x_{2}\})\{x_{3}\}-(x_1\{x_2\}* x_3 -x_3* x_1\{x_2\})\\
&=(x_1* x_2 -x_{2}* x_{1})\{x_{3}\}-x_{1}\{x_{2}\}* x_{3}+x_{3}*
  x_{1}\{x_{2}\}\\
&=x_{1}* x_{2}\{x_{3}\}+x_{1}\{x_{3}\}* x_{2}-
x_{2}* x_{1}\{x_{3}\}\\
&\qquad\qquad\qquad-x_{2}\{x_{3}\}* x_{1}
-x_{1}\{x_{2}\}* x_{3}+x_{3}* x_{1}\{x_{2}\}.
\end{aligned}
\end{equation}
This completes the proof.
\end{proof}

\begin{theorem}\label{t44}
The morphism $I$ is an $S_{n}$--equivariant quasi-isomorphism of complexes, and moreover the following diagram commutes:
\begin{equation}
\xymatrix{
H_{\lf}^{2n-*}(F(n,\R^2))\ar@{=}^\wr[d]\ar@{=}[r]^-{\hbox{\footnotesize$\sim$}}&H_{*}(F(n,\R^2))\ar@{=}^\wr[d]\\
H_{*}(\FF_{2}(n))\ar[r]_{I_{*}}&H_{*}(\SSS_{2}(n)).
}\raisebox{-40pt}[1pt]{}
\label{eq44}
\end{equation}
\end{theorem}

\fullref{t44} will be proven in \fullref{s5}.

A corollary of \fullref{t44} is \fullref{t32}. Indeed, operations $\xi_{1}$, $\zeta_{1}$ correspond to some cycles in
$H_{*}(F(p,\R^2),\pm\Z_{p})$ of degree $p-1$, and $p-2$ respectively. 

Unfortunately there is no proof in the literature 
 of the fact  that $\xi_{1}$ operation of Fred Cohen corresponds exactly to the cycle 
$\epsilon(p)\in\BB_{2}(p,\pm\Z_{p})$ (in the computations of  F\,V~Vainshtein).  At least the author failed to find such a reference. So, we  suppose that
$$
\epsilon(p)=\lambda \xi_1,
$$
for some $\lambda\neq 0$. We will prove in \fullref{s6} that $\lambda=1$.

Modulo the above remark, the cycle $\zeta_{1}=\beta\xi_1$ corresponds to $\beta(\epsilon(p))$. But from~\eqref{eq42'},
$$\beta(\epsilon(p))=\sum_{i=1}^{p-1}\frac{(p-1)!}{i!(p-i)!}\epsilon(i,p-i)=-\sum_{i=1}^{p-1}\frac{(-1)^i}{i}
\epsilon(i,p-i).$$

\section[Proof of Theorem~\ref{t44}]{Proof of \fullref{t44}}\label{s5}

The operads $\Lie$ of Lie algebras, $\prelie$ of pre-Lie algebras, and $\Brace$ of brace algebras, are well known.
A beautiful description of $\prelie$ is given by Chapoton and Livernet~\cite{Cha}. Brace algebras  were introduced by E~Getzler in~\cite{Get}, and by T~Kadeishvili in~\cite{Kad},
see also Gerstenhaber and Voronov~\cite{GV}.
One has  natural inclusions
\begin{equation}
\xymatrix{
\Lie\,\ar@{^{(}->}[r]^-{\iota_{1}}& \prelie\,
\ar@{^{(}->}[r]^-{\iota_{2}}&\Brace,
}
\label{eq51}
\end{equation}
where $[x_{1},x_{2}]$ is mapped to $x_{1}\circ x_{2}-x_{2}\circ x_{1}$; $x_{1}\circ x_{2}$ is mapped to $x_{1}\{x_{2}\}$.

Denote by $\Lie_{1}$,  $\prelie_{1}$, $\Brace_{1}$ the operads of Lie,  pre-Lie,  and Brace algebras with bracket $[\, .\, , \, .\, ]$ of degree one,
pre-Lie product $\circ$ of degree one, brace operations $\{\,\}_{n}$, $n\geq 2$, of degree $n-1$, respectively. One passes from $\Lie$, $\prelie$, $\Brace$ 
algebras to $\Lie_{1}$, $\prelie_{1}$, $\Brace_{1}$ algebras by a desuspension of the underlying spaces.

One has inclusions
 \begin{equation}
\xymatrix{
\Lie_{1}\ar@{^{(}->}[r]^-{\iota_{1}}& \prelie_{1}
\ar@{^{(}->}[r]^-{\iota_{2}}&\Brace_{1}\ar@{^{(}->}[r]^-{\iota_{3}}&\SSS_{2},
}\label{eq52}
\end{equation}
where $\iota_{1}$, $\iota_{2}$ are superanalogues of the inclusions~\eqref{eq51}.

Recall that the homology operad of $\SSS_{2}$ is the operad $\Gerst$ of
Gerstenhaber algebras, ie graded commutative algebras endowed
with a degree one Lie bracket compatible with multiplication
\begin{equation}
[a,bc]=[a,b]c+(-1)^{\deg(b)\cdot \deg(c)}[a,c]b.
\label{eq53}
\end{equation}
It is well known that $\dim\Lie(n)=(n-1)!$, and
$\dim\prelie(n)=n^{n-1}$~\cite{Cha}. Consider an $n!$--dimensional subspace
$\vertical(n)\subset\prelie(n)$ (resp. $\vertical_{1}(n)\subset\prelie_{1}(n)$), which is spanned by the elements
\begin{equation}
(\cdots((x_{\sigma_{1}}\circ x_{\sigma_{2}})\circ x_{\sigma_{3}})\circ\cdots )\circ x_{\sigma_{n}},\label{eq54}
\end{equation}
where $\sigma$ is a permutation from the symmetric group $S_{n}$.

\begin{lemma}\label{l55}
$\iota_{1}(\Lie(n))\subset\vertical(n)$, $\iota_{1}(\Lie_{1}(n))\subset\vertical_{1}(n)$. 
\end{lemma}

\begin{proof}[Proof of \fullref{l55}] To avoid the problem of signs we will consider the first situation.
The second case is obtained by tensoring with the sign representation $\pm \kk$ of the symmetric group.

Decomposition of a bracket from $\Lie(n)$ in the basis~\eqref{eq54} will be called {\it vertical} decomposition.

We will prove our lemma by induction over $n$. For $n=1$, it is evident. Now, suppose $[A,B]\in\Lie(n)$ is some bracket. Then,
\begin{equation}
\iota_{1}([A,B])=\iota_{1}(A)\circ\iota_{1}(B)-\iota_{1}(B)\circ \iota_{1}(A).
\label{eq56}
\end{equation}
We will prove that each summand of~\eqref{eq56} belongs to $\vertical(n)$. To do this, we apply the vertical decomposition for the left
factors of~\eqref{eq56}, and then to each summand of the obtained expression, we apply many times the identity
\begin{equation}
a\circ [b,c]=(a\circ b)\circ c - (-1)^{|b|\cdot |c|}(a\circ c)\circ b,
\label{eq57}
\end{equation}
which is another form of the standard pre-Lie product identity~\cite{Cha}.  For instance, for
$[x_{1},[x_{2},x_{3}]]$, one gets
\begin{multline*}
\iota_{1}([x_{1},[x_{2},x_{3}]])=x_{1}\circ[x_{2},x_{3}]-[x_{2},x_{3}]\circ
x_{1}=\\
(x_{1}\circ x_{2})\circ x_{3}-(x_{1}\circ x_{3})\circ x_{2}-(x_{2}\circ x_{3})\circ x_{1}+ (x_{3}\circ x_{2})\circ x_{1}.
\end{multline*}
This completes the proof of \fullref{l55}.
\end{proof}

Now, we are ready to prove \fullref{t44}. 
Note, that complexes $\FF_{2}(n)$ and $\SSS_{2}(n)$ have the same homology. This homology is
$$
H_*(F(n,\R^2))\simeq \Gerst(n).
$$ 
So, it is sufficient to prove that the induced homology morphism $I_{*}$ is surjective. Consider the maximal degree $n-1$.
In this degree the homology group is $\Lie_{1}(n)$. It follows from \fullref{l55}, that any homology cycle of $\SSS_{2}(n)$ is an image 
of~$I_*$.
For smaller degrees, one needs to use  the cup-product $*$ to obtain all the homology classes as image of~$I_*$. 

To see that diagram~\eqref{eq44} is commutative, one needs to analyze
Poincar\'e duality~\eqref{eq42'''} and the isomorphism $I_{*}$ in
more detail.\hfill\qed

\section{Coefficient $\lambda=1$}\label{s6} 

We have seen in Sections~\ref{s3} and~\ref{s4} that our operation
$\alpha\co x\mapsto x^{[p]}$ 
defined on odd degree cycles of a Hochschild complex (and on cycles of any degree if $p=2$) 
is a multiple of  Cohen's $\xi_1$:
\begin{equation}
\alpha(x)=\lambda\xi_1(x)
\label{eq60}
\end{equation}
We need to prove that $\lambda$ is in fact one.

Notice that $\alpha=\lambda\xi_1$ is a cycle of maximal degree $*=p-1$ in the complex $\SSS_2(p)\otimes_{S_p}(\pm\Z_p)$. The same is true for the Gerstenhaber bracket $[\, .\, ,\, . \, ]\in (\SSS_2(2))_1$.  It means that each of these two cycles  has only one representative in the corresponding complex. Any composition of these two operations is also a cycle
of maximal degree in some complex
$\SSS_{2}(n)\otimes_{S_{k_1}\times\cdots\times S_{k_\ell}}W$, where $W$ is of type~\eqref{eq04}, and hence also has a unique representative.

In~\cite{Co} F~Cohen proved the following identity (for the operations
induced by $\CC_2$--action):
\begin{equation}
[y,\xi_1(x)]=[\ldots[[y,\underbrace{x],x]\ldots x]}_{p}
\label{eq61}
\end{equation}
Therefore one has the identity
\begin{equation}
[y,x^{[p]}]=\lambda [\ldots[[y,\underbrace{x],x]\ldots x]}_{p},
\label{eq62}
\end{equation}
for some $\lambda\neq 0$. But note that identity~\eqref{eq62} is on the level of chains. Indeed, by~\eqref{eq60}, and~\eqref{eq61}
the left-hand side and the right-hand side of~\eqref{eq62} are representatives of the same cycle in the homology of the complex
$\SSS_2(p+1)\otimes_{S_p\times S_1}\sign_1$, where $\sign_1$ is the sign
representation of $S_p$. But this homology class lies in the maximal
degree $\deg=p$ of the above complex, so there exists only one its representative.

As a consequence~\eqref{eq62} is true 
for any elements $x$ and $y$ of any $\SSS_2$--algebra (we assume  that
$p(\deg(x)-1)$ is even).
Hence it is true for any $\Brace_1$--algebra, and hence for any $\Brace$--algebra (except that $p\cdot
\deg(x)$ is now even). 

Because of the natural map
$$
\Brace\longrightarrow\Assoc
$$
(the operation $\{\,\}_2$ is mapped to the product, all the other braces $\{\,\}_n$, $n\geq 3$, are mapped to zero), any associative algebra can
be considered as a brace algebra. For associative algebras $x^{[k]}=x^k$, and
the identity~\eqref{eq62} is known to be true with the coefficient $\lambda=1$, cf \cite[Chapter~V]{J}. 
We obtain as a consequence that $\lambda=1$.

\section{Some relations that hold in maximal degree}\label{s7}
In this section we use grading $|\, .\,|=\deg+1$. So, Lie, pre-Lie and brace operations become of degree zero.

Arguing as in the previous section we can discover many other interesting identities that hold already on the level of chains.
For instance, the Gerstenhaber bracket~\eqref{eq75} satisfies the
Jacobi identity
\begin{equation}
(-1)^{|a_{1}|\cdot|a_{3}|}[[a_{1},a_{2}],a_{3}]+(-1)^{|a_{2}|\cdot|a_{1}|}[[a_{2},a_{3}],a_{1}]+(-1)^{|a_{3}|\cdot|a_{2}|}[[a_{3},a_{1}],a_{2}],
\label{eq70}
\end{equation}
and the identities
\begin{align}
[x,x]&=0 \quad\text{for $p=2$,}\label{eq71}\\
[[x,x],x]&=0 \quad\text{for $p=3$.}
\label{eq72}
\end{align}
In~\cite{T5} the author proved the identities
\begin{gather}
(a^{[p^k]})^{[p^l]}=a^{[p^{k+l}]},\label{eq73}\\
(b^{[2p^k]})^{[p^l]}=b^{[2p^{k+l}]},\label{eq74}
\end{gather}
whenever $|a|$ is even, and $|b|$ is odd. The proof was a tedious check  of some combinatorial properties of planar rooted trees.

All the above relations hold for any brace algebra. The following theorem provides more other relations.

\begin{theorem}\label{t75}
Any brace algebra in characteristic $p$ is  a $p$--restricted Lie algebra with bracket
\begin{equation}
[a,b]:=a\{ b\}-(-1)^{|a|\cdot |b|}b\{a\},
\label{eq75}
\end{equation}
and restriction operation (defined for elements of even degree $|\, . \,|$)
$$
a^{[p]}:=a\underbrace{\{a\}\ldots\{a\}}_{p-1}. 
$$
\end{theorem}

The above theorem means that Jacobson's relations hold (see
Jacobson~\cite[Section~V.7]{J}), that is,
\begin{gather}
[a,b^{[p]}]=[\ldots[[a,\underbrace{b],b]\ldots b]}_{p},
\label{eq76}\\
(c_{1}+c_{0})^{[p]}=c_{1}^{[p]}+c_{0}^{[p]}+\sum_{i=1}^{p-1}d_{i}(c_{1},c_{0}),
\label{eq77}
\end{gather}
where 
$$
i\cdot d_{i}(c_{1},c_{0})=\sum_{ {\epsilon_{i}\in\{0,1\}, \, i=1\dots p-2} \atop {\epsilon_{1}+\dots \epsilon_{p-2}=i-1} }
[\ldots[[[c_{1},c_{0}],c_{\epsilon_{1}}],c_{\epsilon_{2}}],\ldots,c_{\epsilon_{p-2}}].
$$
The elements $b$, $c_{1}$, $c_{0}$ are even.

All the above relations arise as a manifestation of the fact that homology classes in maximal degree of the complexes
$$
\SSS_{2}(n)\otimes_{S_{k_1}\times S_{k_2}\times\cdots\times S_{k_\ell}} W
$$
have unique representatives. Composition of operations in maximal degree is also an operation in maximal degree.
So, the relations \eqref{eq70}, \eqref{eq71}, \eqref{eq72}, \eqref{eq73}, \eqref{eq74},  \eqref{eq76}, \eqref{eq77} 
follow from the analogous relations
for homology operations of double loop spaces (see Cohen \cite{Co}), and  from
\fullref{c21}.

Now, note that to define the above operations we need only pre-Lie product $\circ$. We define 
$$
x^{[k]}:=(\dots((\underbrace{x\circ x)\circ x)\dots)\circ x}_{\text{$k$ times}}.
$$ 
So, it is natural to ask
whether these relations hold for pre-Lie algebras. (For  brace algebras all the relations hold automatically.)

\begin{theorem}\label{t78}
{\rm (a)}\qua Relations \eqref{eq70}, \eqref{eq71}, \eqref{eq72} hold for any pre-Lie algebra.

{\rm (b)}\qua Relations \eqref{eq73}, \eqref{eq74} do not hold for a free pre-Lie algebra with one even generator $a$,
resp. one odd generator  $b$.

{\rm (c)}\qua Relation \eqref{eq76} does not hold for a free pre-Lie algebra with two generators $a$ and $b$ {\rm (}the second one being even{\rm )}.

{\rm (d)}\qua Relation \eqref{eq77} holds for any pre-Lie algebra. 
\end{theorem}

\begin{proof} (a) is well known. (b) and (c) are easy to verify if one uses the representation of free pre-Lie algebras
in terms of rooted trees~\cite{Cha}. In fact, equality~\eqref{eq76}  holds if one adds to the left hand-side the  rooted tree

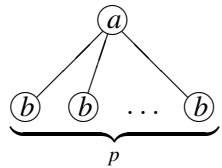
\begin{center}
\unitlength=0.2em
 \begin{picture}(30,22)(0,-5)
   \put(0,0){\circle{5}}
   \put(2.0,1.7){\line(1,1){11.3}}
   \put(10,0){\circle{5}}
   \put(10.79,2.37){\line(1,3){3.42}}
   \put(30,0){\circle{5}}
   \put(28,1.7){\line(-1,1){11.3}}
   \put(15,15){\circle{5}}
   \put(-1.30,-1.60){$b$}
   \put(8.7,-1.6){$b$}
   \put(17,-1.2){$\ldots$}
   \put(28.7,-1.6){$b$}
   \put(13.7,13.6){$a$}
   \put(-3,-3){$\underbrace{\hspace{7.2em}}_{p}$}
    
 \end{picture}
\end{center}

The proof of (d) is a slight modification of the proof of the same result for associative algebras given in~\cite[Section~V.7]{J}.
\end{proof}

\begin{lemma}\label{l79}
For any  even elements $a$, $b$ of a pre-Lie algebra the following identities hold:
$$
[\ldots[[a,\underbrace{b],b],\ldots b]}_{N}=\sum_{i=0}^{N}(-1)^{i}{N\choose{i}}(\cdots(\underbrace{b\circ b)\circ b)\cdots
\circ b}_{i})\circ a)\circ
\underbrace{b)\circ b)\cdots )\circ b}_{N-i-1}. 
$$
\end{lemma}

\begin{proof}
Induction over $N$.
\end{proof}

In the case of characteristic $p$ and $N=p-1$, one gets:
\begin{equation}
[\ldots[[a,\underbrace{b],b],\ldots b]}_{p-1}=\sum_{i=0}^{p-1}(\cdots(\underbrace{b\circ b)\circ b)\cdots
\circ b}_{i})\circ a)\circ
\underbrace{b)\circ b)\cdots )\circ b}_{p-i-1}.
\label{eq710}
\end{equation}

Consider equality~\eqref{eq710} for $a=c_{1}$, $b=c_{1}+\lambda c_{0}$, and differentiate it over $\lambda$. Identity \eqref{eq77}~follows from the obtained
expression.

\begin{remark}[M~Livernet]
\label{r64}
Actually \fullref{t75} is an obvious consequence of a result of
B~Fresse~\cite{Fresse1,Fresse2}. The usual definition of an algebra over
an operad $\O$
consists in defining a family of compatible maps
$$
(\O(n)\otimes V^{\otimes n})_{S_{n}}\to V,
$$
where $(-)_{S_{n}}$ denotes the space of coinvariants of the symmetric
group action. Instead of doing this, B~Fresse proposes to define
$\Gamma\O$--action as 
a family of compatible maps
$$
(\O(n)\otimes V^{\otimes n})^{S_{n}}\to V,
$$
where $(-)^{S_{n}}$ is now the space of invariants of the symmetric group action. In  case $\kk$ is not a field of characteristic zero, this defines
another algebraic structure, which Fresse calls an $\O$--algebra with divided symmetries. For example, if $\O=\Comm$ is the operad of commutative algebras, then a
$\Gamma\Comm$--algebra is a so called divided system -- a commutative algebra with unary divided power operations. 
Fresse proves that for $char(\kk)=p$, the $\Gamma\Lie$--algebra structure is exactly the restricted Lie algebra structure. 

On the other hand, the $S_{n}$--action on $\Brace(n)$ is free. Hence,
$\Brace$ and $\Gamma\Brace$ are the same structures. But any
$\Gamma\Brace$--algebra
must be $\Gamma\Lie$ due to the morphism of operads:
$$
\Lie\to\Brace.
$$
This proves the theorem. 

Note, that this argument does not work for the operad $\prelie$, since
$\prelie(n)$ are not projective $S_{n}$--modules for $n\geq 3$.
\end{remark}

The following proposition provides one more relation which holds on the level of chains.

\begin{proposition}\label{p65}
Consider a Hochschild complex obtained from another Hochschild complex (free over $\Z$) by tensoring with $\Z_p$, $p\geq 3$. Suppose $x$ is a representative of an odd degree cycle in this complex. The following
identity holds on the level of chains:
$$
\zeta_{1}x=\beta(\xi_{1}x)- (\ad^{p-1}x)(\beta x),
$$
where operator $\ad\, x$ is the adjoint action $[x,\, .\, ]$.
\end{proposition}

\begin{proof} It is a consequence of~\eqref{eq710} applied for $a=\beta x$, and $b=x$. 
\end{proof} 

\section{About higher Dyer--Lashof--Cohen operations and Steenrod powers}\label{s8}
Homology operations for iterated $d$--loop spaces correspond to
equivariant cycles of the operad $\CC_{d}$ of little $d$--cubes.
Following F~Cohen~\cite{Co}, these operations are generated by the Pontryagin product $*$, the higher Browder operator $[\, .\, ,\, .\, ]_{d-1}$,
and operations $\xi_{i}$, $\zeta_{i}$, $i=1\ldots d-1$:
$$
\begin{array}{ccl}
* &\in&H_{0}(\CC_{d}(2))=H_{0}(S^{d-1})\\
{[\, .\, , \, . \, ]}&\in&H_{d-1}(\CC_{d}(2))=H_{d-1}(S^{d-1})\\
\xi_{i}&\in&H_{i(p-1)}(\CC_{d}(p)/S_{p},(\pm\Z_{p})^{\otimes i})\\
\zeta_{i}=\beta(\xi_{i})&\in&H_{i(p-1)-1}(\CC_{d}(p)/S_{p},(\pm\Z_{p})^{\otimes i}),\,  p>2
\end{array}
$$
Operations $*$; $\xi_{i}$, $\zeta_{i}$, $i=1\ldots d-2$, are inherited
from the $(d-1)$--cubes action. The Browder bracket $[\, .\, ,\, .\, ]_{d-2}$ defined for
$(d-1)$--loops becomes trivial for $d$--loops.

While $d=\infty$, the Browder operator disappears and we have only $*$; $\xi_{i}$, $\zeta_{i}$, $i\in\N$.

J\,E~McClure and J\,H~Smith defined a differential graded operad $\SSS$
(\emph{operad of surjections}), which is naturally filtered,
\begin{equation}
\SSS_{1}\subset\SSS_{2}\subset\SSS_{3}\subset\cdots=\SSS,
\label{eqFiltr}
\end{equation}
and whose $d^{\text{th}}$ filtration term $\SSS_{d}$ is quasi-isomorphic to the
singular chains operad of little $d$--cubes~\cite{McSm2}.\footnote{See also~\cite{BF,Smith}.}
The operad $\SSS_{2}$, considered in this paper, acts on Hochschild complexes. The whole operad $\SSS$ acts on singular cochains of topological spaces.
Operations $\xi_{i}$, $i\in\N$, in the last situation, are the well known Steenrod powers.

Homology operations induced by an $\SSS_{d}$--action  correspond to  equivariant cycles
of $\SSS_{d}(n)$,
$$
H_{*}(\SSS_{d}(n)\otimes_{S_{k_1}\times\cdots\times S_{k_\ell}}W)\simeq H_{*}(\CC_{d}(n)/S_{k_1}\times\cdots\times S_{k_\ell},W),
$$
where $W$ is as usual of the form~\eqref{eq04}.

Cycles $*$, $[\, .\, ,\, .\, ]_{d-1}$ are  defined as
\begin{equation}
\begin{aligned}
*&=x_{1}\cup_{0}x_{2},\\
[\, .\, ,\, .\, ]_{d-1}&=x_{1}\cup_{d-1}x_{2}-x_{2}\cup_{d-1}x_{1},\label{eqBr}
\end{aligned}
\end{equation}
where $\cup_{d-1}$ denotes the $(d-1)$--cup product, see Steenrod
\cite{Steen}, and McClure--Smith \cite{McSm2}.

In the case $p=2$, all the operations $\xi_{i}$, $i=1,2,3,\ldots$ are also
explicit. They were defined by N\,E~Steenrod in his seminal work~\cite{Steen}:
$$
\xi_{i}(x)=x\cup_{i}x
$$
See also~\cite{McSm2} for a description of the cup-products $\cup_{i}\in\SSS_{i+1}(2)$ as elements of the operad of surjections.

In the case $p>2$, all the constructions are implicit.\footnote{Except the
construction of Gonz\'alez-D\'{\i}az and Real \cite{GR}, that unfortunately does not respect the filtration~\eqref{eqFiltr}.} One knew only $\xi_{0}(x)$ which is $x^p$:
$$
\xi_{0}(x)=\underbrace{x\cup_{0}x\cup_{0}\ldots\cup_{0}x}_{\text{$p$ times}}
$$
($\cup_{0}$--product is associative). The result of this paper permits to define explicitly the 
operation~$\xi_{1}$:
$$
\xi_{1}(x)=(\dots((\underbrace{x\cup_{1} x)\cup_{1} x)\dots)\cup_{1} x}_{\text{$p$ times}},
$$
which is the last non-trivial Steenrod power defined for odd degree elements:
$$
\xi_{1}(x)=P^{\frac{\deg x-1}2}x,
$$
all the higher Steenrod powers $P^{i}$, $i>\frac{\deg x-1}2$, are trivial:
$P^{i}x=0$ (see Steenrod and Epstein \cite{SteenEp}).

The question is whether it is possible to find explicit formulae for all the other $\xi_{i}$, $i\geq 2$, $p\geq 3$.
Note that the main difficulty is to find these formulae, since operations $\zeta_{i}$, $i\geq 1$, as cycles of $\SSS(p)\otimes_{S_{p}}(\pm\Z_{p})^{\otimes i}$,
are Bockstein images  $\beta(\xi_{i})$, $i\geq 1$.

Unfortunately the methods given in this paper can not be generalized in an obvious way to define explicitly $\xi_{i}$, $i\geq 2$, $p\geq 3$.

Of course, one can easily define complexes $\FF_{d}(n)$, for any $d\geq 1$ and $n\geq 1$, taking the so called \lq\lq lexicographical cellular decomposition" of
the one-point compactification of the configuration spaces $F(n,\R^d)$. Operation $\xi_{d-1}$ corresponds to the cycle spanned by the only cell
of $\overline{B(p,\R^d)}=\overline{F(p,\R^d)/S_{p}}$: all  $p$ points having the same first $d-1$ coordinates. But the author did not manage to define
an analogous quasi-isomorphism $\FF_{d}\to\SSS_{d}$ for $d\geq 3$.

The results of Fresse, see \fullref{r64},  also can not be applied: higher Browder operators~\eqref{eqBr}, $d\geq 3$, satisfy
the Jacobi identity only up to a non-trivial homotopy. Hence we can not define a map $\Lie_{d-1}\to\SSS_{d}$, where $\Lie_{d-1}$ denotes the 
operad of Lie algebras with the bracket of degree $d-1$, the bracket being sent to the Browder operator.

\bibliographystyle{gtart}
\bibliography{link}

\begin{thebibliography}{}
\providecommand\bibmarginpar{\leavevmode\marginpar}
\def\urlstyle#1{{\tt #1}}

\bibitem{BF}
\textbf{C Berger}, \textbf{B Fresse},
  \href{http://dx.doi.org/10.1017/S0305004103007138} {\emph{Combinatorial
  operad actions on cochains}}, Math. Proc. Cambridge Philos. Soc. 137 (2004)
  135--174 \xox{MR}{2075046}

\bibitem{Cha}
\textbf{F Chapoton}, \textbf{M Livernet},
  \href{http://dx.doi.org/10.1155/S1073792801000198} {\emph{Pre-{L}ie algebras
  and the rooted trees operad}}, Internat. Math. Res. Notices  (2001) 395--408
  \xox{MR}{1827084}

\bibitem{Co}
\textbf{F\,R Cohen}, \emph{The homology of $\mathcal{C}_{n+1}$--spaces, $n \geq
  0$}, from: ``The homology of iterated loop spaces'', Lecture Notes in
  Mathematics 533, Springer, Berlin (1976)  207--351 \xox{MR}{0436146}

\bibitem{DyLash}
\textbf{E Dyer}, \textbf{R\,K Lashof}, \emph{Homology of iterated loop spaces},
  Amer. J. Math. 84 (1962) 35--88 \xox{MR}{0141112}

\bibitem{Fresse1}
\textbf{B Fresse}, \href{http://dx.doi.org/10.1016/S0764-4442(97)83949-3}
  {\emph{Op\'erations de {C}artan pour les alg\`ebres simpliciales sur une
  op\'erade}}, C. R. Acad. Sci. Paris S\'er. I Math. 325 (1997) 247--252
  \xox{MR}{1464814}

\bibitem{Fresse2}
\textbf{B Fresse}, \href{http://dx.doi.org/10.1090/S0002-9947-99-02489-7}
  {\emph{On the homotopy of simplicial algebras over an operad}}, Trans. Amer.
  Math. Soc. 352 (2000) 4113--4141 \xox{MR}{1665330}

\bibitem{F}
\textbf{D\,B Fuchs}, \emph{Cohomology of braid groups ${\rm mod}2$}, Funct.
  Analysis and Appl. 4 (1970) 62--73

\bibitem{Ge}
\textbf{M Gerstenhaber}, \emph{The cohomology structure of an associative
  ring}, Ann. of Math. $(2)$ 78 (1963) 267--288 \xox{MR}{0161898}

\bibitem{GV}
\textbf{M Gerstenhaber}, \textbf{A\,A Voronov},
  \href{http://dx.doi.org/10.1155/S1073792895000110} {\emph{Homotopy
  {$G$}--algebras and moduli space operad}}, Internat. Math. Res. Notices
  (1995) 141--153 \xox{MR}{1321701}

\bibitem{Get}
\textbf{E Getzler}, \emph{Cartan homotopy formulas and the {G}auss-{M}anin
  connection in cyclic homology}, from: ``Quantum deformations of algebras and
  their representations (Ramat-Gan, 1991/1992; Rehovot, 1991/1992)'', Israel
  Math. Conf. Proc. 7, Bar-Ilan Univ., Ramat Gan (1993)  65--78
  \xox{MR}{1261901}

\bibitem{GR}
\textbf{R Gonz{\'a}lez-D{\'{\i}}az}, \textbf{P Real},
  \href{http://dx.doi.org/10.1016/S0022-4049(99)00006-7} {\emph{A combinatorial
  method for computing {S}teenrod squares}}, J. Pure Appl. Algebra 139 (1999)
  89--108 \xox{MR}{1700539}

\bibitem{J}
\textbf{N Jacobson}, \emph{Lie algebras}, Interscience Tracts in Pure and
  Applied Mathematics, No. 10, Interscience Publishers (a division of John
  Wiley \& Sons), New York-London (1962) \xox{MR}{0143793}

\bibitem{Kad}
\textbf{T\,V Kadeishvili}, \emph{The structure of the {$A(\infty)$}--algebra,
  and the {H}ochschild and {H}arrison cohomologies}, Trudy Tbiliss. Mat. Inst.
  Razmadze Akad. Nauk Gruzin. SSR 91 (1988) 19--27 \xox{MR}{1029003}

\bibitem{K2}
\textbf{M Kontsevich}, \href{http://dx.doi.org/10.1023/A:1007555725247}
  {\emph{Operads and motives in deformation quantization}}, Lett. Math. Phys.
  48 (1999) 35--72 \xox{MR}{1718044}

\bibitem{KS}
\textbf{M Kontsevich}, \textbf{Y Soibelman}, \emph{Deformations of algebras
  over operads and the {D}eligne conjecture}, from: ``Conf\'erence Mosh\'e
  Flato 1999, Vol. I (Dijon)'', Math. Phys. Stud. 21, Kluwer Acad. Publ.,
  Dordrecht (2000)  255--307 \xox{MR}{1805894}

\bibitem{AK2}
\textbf{T Kudo}, \textbf{S Araki}, \emph{Topology of $H_n$--spaces and
  $H$--squaring operations}, Mem. Fac. Sci. Ky\=usy\=u Univ. Ser. A. 10 (1956)
  85--120 \xox{MR}{0087948}

\bibitem{Ma}
\textbf{N\,S Markaryan}, \href{http://dx.doi.org/10.1007/BF02307210}
  {\emph{Homology of braid groups with nontrivial coefficients}}, Mat. Zametki
  59 (1996) 846--854, 960 \xox{MR}{1445470}

\bibitem{McSm1}
\textbf{J\,E McClure}, \textbf{J\,H Smith}, \emph{A solution of {D}eligne's
  {H}ochschild cohomology conjecture}, from: ``Recent progress in homotopy
  theory (Baltimore, MD, 2000)'', Contemp. Math. 293, Amer. Math. Soc.,
  Providence, RI (2002)  153--193 \xox{MR}{1890736}

\bibitem{McSm2}
\textbf{J\,E McClure}, \textbf{J\,H Smith},
  \href{http://dx.doi.org/10.1090/S0894-0347-03-00419-3} {\emph{Multivariable
  cochain operations and little {$n$}--cubes}}, J. Amer. Math. Soc. 16 (2003)
  681--704 \xox{MR}{1969208}

\bibitem{Smith}
\textbf{J\,H Smith}, \emph{Simplicial group models for {$\Omega\sp nS\sp nX$}},
  Israel J. Math. 66 (1989) 330--350 \xox{MR}{1017171}

\bibitem{Steen}
\textbf{N\,E Steenrod},
  \href{http://links.jstor.org/sici?sici=0003-486X(194704)2:48:2%3C290:POCAEO%%
3E2.0.CO%3B2-A} {\emph{Products of cocycles and extensions of mappings}}, Ann.
  of Math. $(2)$ 48 (1947) 290--320 \xox{MR}{0022071}

\bibitem{SteenEp}
\textbf{N\,E Steenrod}, \textbf{D\,B\,A Epstein}, \emph{Cohomology operations},
  Annals of Mathematics Studies 50, Princeton University Press, Princeton, N.J.
  (1962) \xox{MR}{0145525}

\bibitem{Tam1}
\textbf{D\,E Tamarkin}, \emph{Another proof of M Kontsevich formality theorem}
  \xox{arXiv}{math.QA/9803025}

\bibitem{Tam2}
\textbf{D\,E Tamarkin},
  \href{http://dx.doi.org/10.1023/B:MATH.0000017651.12703.a1} {\emph{Formality
  of chain operad of little discs}}, Lett. Math. Phys. 66 (2003) 65--72
  \xox{MR}{2064592}

\bibitem{T5}
\textbf{V Tourtchine}, \emph{On the other side of the bialgebra of chord
  diagrams}, J. Knot Th. Ram. (to appear) \xox{arXiv}{math.QA/0411436}

\bibitem{Vain}
\textbf{F\,V Va{\u\i}n{\v{s}}te{\u\i}n}, \emph{The cohomology of braid groups},
  Funktsional. Anal. i Prilozhen. 12 (1978) 72--73 \xox{MR}{498903}

\bibitem{V0}
\textbf{V\,A Vassiliev}, \href{http://dx.doi.org/10.1007/BF01077624}
  {\emph{Cohomology of braid groups and the complexity of algorithms}},
  Funktsional. Anal. i Prilozhen. 22 (1988) 15--24, 96 \xox{MR}{961758}

\bibitem{Vor}
\textbf{A\,A Voronov}, \emph{Homotopy {G}erstenhaber algebras}, from:
  ``Conf\'erence Mosh\'e Flato 1999, Vol. II (Dijon)'', Math. Phys. Stud. 22,
  Kluwer Acad. Publ., Dordrecht (2000)  307--331 \xox{MR}{1805923}

\bibitem{Weibel}
\textbf{C\,A Weibel}, \emph{An introduction to homological algebra}, Cambridge
  Studies in Advanced Mathematics 38, Cambridge University Press, Cambridge
  (1994) \xox{MR}{1269324}

\bibitem{Craig}
\textbf{C Westerland}, \href{http://dx.doi.org/10.1007/s00209-005-0778-9}
  {\emph{Dyer--Lashof operations in the string topology of spheres and
  projective spaces}}, Math. Z. 250 (2005) 711--727 \xox{MR}{2179618}

\end{thebibliography}

\end{document}